\documentclass[12pt,twoside]{article}

\usepackage[english]{babel}
\usepackage{umj_1}           

\usepackage{graphicx, srcltx}   
\usepackage{cite}               

\usepackage{amsmath}            
\usepackage{amsfonts,amssymb}   
\usepackage{srcltx}             

\newcommand{\Leftclt}{Roman A. Veprintsev}
\newcommand{\Rightclt}{On Hardy\,--\,Littlewood-type and Hausdorff\,--\,Young-type inequalities}
\usepackage{authblk}

\DeclareMathOperator*{\esssup}{ess\,sup}

\begin{document}

\begin{Titul}
{\large \bf ON HARDY\,--\,LITTLEWOOD-TYPE\\ AND HAUSDORFF\,--\,YOUNG-TYPE INEQUALITIES\\[0.2em] FOR GENERALIZED GEGENBAUER EXPANSIONS }\\[3ex]
{{\bf Roman~A.~Veprintsev} \\[5ex]}
\end{Titul}

\begin{Anot}
{\bf Abstract.} Using well-known techniques, we establish Hardy\,--\,Littlewood-type and Hausdorff\,--\,Young-type inequalities for generalized Gegenbauer expansions and their unification.


\end{Anot}


\section{Introduction and preliminaries}

In this section, we introduce some classes of orthogonal polynomials on $[-1,1]$, including the so-called generalized Gegenbauer polynomials. For a background and more details on the orthogonal polynomials, the reader is referred to \cite{dai_xu_book_approximation_theory_2013,dunkl_xu_book_orthogonal_polynomials_2014,andrews_askey_roy_book_special_functions_1999,szego_book_orthogonal_polynomials_1975}.

Let $\alpha,\,\beta>-1$. The Jacobi polynomials, denoted by $P_n^{(\alpha,\beta)}(\cdot)$, where $n=0,1,\ldots$, are orthogonal with respect to the Jacobi weight function $w_{\alpha,\beta}(t)=(1-t)^\alpha(1+t)^\beta$ on $[-1,1]$, namely,
\begin{equation*}
\int\nolimits_{-1}^1 P_n^{(\alpha,\beta)}(t)\, P_m^{(\alpha,\beta)}(t)\, w_{\alpha,\beta}(t)\,dt=\begin{cases}\dfrac{2^{\alpha+\beta+1}\Gamma(n+\alpha+1)\Gamma(n+\beta+1)}{(2n+\alpha+\beta+1)\Gamma(n+1)\Gamma(n+\alpha+\beta+1)},&n=m,\\
0,&n\not=m.
\end{cases}
\end{equation*}
Here, as usual, $\Gamma$ is the gamma function.

For $\lambda>-\frac{1}{2}$, $\mu\geq0$, and $n=0,1,\ldots$, the generalized Gegenbauer polynomials $C_n^{(\lambda,\mu)}(\cdot)$ are defined by
\begin{equation*}\label{coefficients_for_generalized_Gegenbauer_polynomials}
\begin{array}{ll}
C_{2n}^{(\lambda,\mu)}(t)=a_{2n}^{(\lambda,\mu)}P_n^{(\lambda-1/2,\mu-1/2)}(2t^2-1), & a_{2n}^{(\lambda,\mu)}=\dfrac{(\lambda+\mu)_n}{(\mu+\frac{1}{2})_n},\\[1.0em]
C_{2n+1}^{(\lambda,\mu)}(t)=a_{2n+1}^{(\lambda,\mu)}\,t P_n^{(\lambda-1/2,\mu+1/2)}(2t^2-1),\quad & a_{2n+1}^{(\lambda,\mu)}=\dfrac{(\lambda+\mu)_{n+1}}{(\mu+\frac{1}{2})_{n+1}},
\end{array}
\end{equation*}
where $(\lambda)_n$ denotes the Pochhammer symbol given by
\begin{equation*}
(\lambda)_0=1,\quad (\lambda)_n=\lambda(\lambda+1)\cdots(\lambda+n-1)\quad\text{ for}\quad n=1,2,\ldots.
\end{equation*}
They are orthogonal with respect to the weight function
\begin{equation}\label{weight_function}
v_{\lambda,\mu}(t)=|t|^{2\mu}(1-t^2)^{\lambda-1/2},\quad t\in[-1,1].
\end{equation}

For $\mu=0$, these polynomials, denoted by $C_n^{\lambda}(\cdot)$, are called the Gegenbauer polynomials:
\begin{equation*}
C_n^{\lambda}(t)=C_n^{(\lambda,0)}(t)=\frac{(2\lambda)_n}{(\lambda+\frac{1}{2})_n} P_n^{(\lambda-1/2,\lambda-1/2)}(t).
\end{equation*}

For $\lambda>-\frac{1}{2}$, $\mu>0$, and $n=0,1,\ldots$, we have the following connection:
\begin{equation*}
C_n^{(\lambda,\mu)}(t)=c_\mu\int\nolimits_{-1}^1 C_n^{\lambda+\mu}(tx)(1+x)(1-x^2)^{\mu-1}\,dx,\quad c_\mu^{-1}=2\int\nolimits_0^1 (1-x^2)^{\mu-1}\,dx.
\end{equation*}

Denote by $\bigl\{\widetilde{C}_n^{(\lambda,\mu)}(\cdot)\bigr\}_{n=0}^{\infty}$ the sequence of orthonormal generalized Gegenbauer polynomials. It is easily verified that these polynomials are given by the following formulae:
\begin{equation*}\label{coefficients_for_orthonormal_generalized_Gegenbauer_polynomials}
\begin{split}
&\widetilde{C}_{2n}^{(\lambda,\mu)}(t)=\widetilde{a}_{2n}^{\,(\lambda,\mu)}P_n^{(\lambda-1/2,\mu-1/2)}(2t^2-1),\\ &\widetilde{a}_{2n}^{\,(\lambda,\mu)}=\Biggl(\dfrac{(2n+\lambda+\mu)\Gamma(n+1)\Gamma(n+\lambda+\mu)}{\Gamma(n+\lambda+\frac{1}{2})\Gamma(n+\mu+\frac{1}{2})}\Biggr)^{1/2},\\[0.4em]
&\widetilde{C}_{2n+1}^{(\lambda,\mu)}(t)=\widetilde{a}_{2n+1}^{\,(\lambda,\mu)}\,t P_n^{(\lambda-1/2,\mu+1/2)}(2t^2-1),\\ &\widetilde{a}_{2n+1}^{\,(\lambda,\mu)}=\Biggl(\dfrac{(2n+\lambda+\mu+1)\Gamma(n+1)\Gamma(n+\lambda+\mu+1)}{\Gamma(n+\lambda+\frac{1}{2})\Gamma(n+\mu+\frac{3}{2})}\Biggr)^{1/2}.
\end{split}
\end{equation*}

The generalized Gegenbauer polynomials play an important role in Dunkl harmonic analysis (see, for example, \cite{dunkl_xu_book_orthogonal_polynomials_2014,dai_xu_book_approximation_theory_2013}). So, the study of these polynomials and their applications is very natural.

The notation $f(n)\asymp g(n)$, $n\to\infty$, means that there exist positive constants $C_1$, $C_2$, and a positive integer $n_0$ such that $0\leq C_1 g(n)\leq f(n)\leq C_2 g(n)$ for all $n\geq n_0$.
For brevity, we will omit ``$n\to\infty$'' in the asymptotic notation.

Define the uniform norm of a continuous function $f$ on $[-1,1]$ by
\begin{equation*}
\|f\|_{\infty}=\max\limits_{-1\leq t\leq 1} |f(t)|.
\end{equation*}
The maximum of two real numbers $x$ and $y$ is denoted by $\max(x,y)$.

In \cite{veprintsev_preprint_max_value_2015}, we prove the following result.

\begin{teoen}\label{main_result_theorem_of_max_value_preprint_for_orthonormal_system}
Let $\lambda>-\frac{1}{2}$, $\mu>0$. Then
\begin{equation*}
\bigl\|\widetilde{C}_n^{(\lambda,\mu)}\bigr\|_{\infty}\asymp n^{\max(\lambda,\mu)}.
\end{equation*}
\end{teoen}

Given $1\leq p\leq\infty$, we denote by $L_p(v_{\lambda,\mu})$ the space of complex-valued Lebesgue measurable functions $f$ on $[-1,1]$ with finite norm
\begin{equation*}
\begin{array}{lr}
\|f\|_{L_p(v_{\lambda,\mu})}=\Bigl(\int\nolimits_{-1}^1 |f(t)|^p\,v_{\lambda,\mu}(t)\,dt\Bigr)^{1/p},&\quad 1\leq p<\infty,\\[1.0em]
\|f\|_{L_\infty}=\esssup\limits_{x\in[-1,1]} |f(x)|,& p=\infty.
\end{array}
\end{equation*}
For a function $f\in L_p(v_{\lambda,\mu})$, $1\leq p\leq\infty$, the generalized Gegenbauer expansion is defined by
\begin{equation*}
f(t)\sim\sum\limits_{n=0}^\infty \hat{f}_n \widetilde{C}_n^{(\lambda,\mu)}(t),\qquad \text{where}\quad\hat{f}_n=\int\nolimits_{-1}^1 f(t)\, \widetilde{C}_n^{(\lambda,\mu)}(t)\,v_{\lambda,\mu}(t)\,dt.
\end{equation*}

For $1<p<\infty$, we denote by $p'$ the conjugate exponent to $p$, that is, $\frac{1}{p}+\frac{1}{p'}=1$.

The aim of this paper is to establish Hardy\,--\,Littlewood-type and Hausdorff\,--\,Young-type inequalities for generalized Gegenbauer expansions in Sections \ref{section_for_Hardy-Littlewood_inequality} and \ref{section_for_Hausdorff-Young_inequality}, respectively. Also, we give their unification in Section \ref{section_for_unification_of_inequalities}.

\section{Hardy\,--\,Littlewood-type inequalities\\ for generalized Gegenbauer expansions}\label{section_for_Hardy-Littlewood_inequality}

The analogue of the Hardy\,--\,Littlewood inequality is given in the following theorem, which can be deduced as a corollary from \cite[Theorems 3.2 and 3.6]{stein_weiss_article_interpolation_1958} (for \eqref{first_part_of_Hardy-Littlewood_inequality} and \eqref{second_part_of_Hardy-Littlewood_inequality}, respectively). Nevertheless, for convenience we give a direct proof of the theorem, based on Theorem \ref{main_result_theorem_of_max_value_preprint_for_orthonormal_system} and our settings.

\begin{teoen}\label{Hardy-Littlewood_inequality_for_generalized_Gegenbauer_expansions}
$(a)$ If $1<p\leq2$ and $f\in L_p(v_{\lambda,\mu})$, then
\begin{equation}\label{first_part_of_Hardy-Littlewood_inequality}
\Bigl\{\sum\limits_{n=0}^\infty\Bigl((n+1)^{\left(\frac{1}{p'}-\frac{1}{p}\right)\left(\max(\lambda,\mu)+1\right)}|\hat{f}_n|\Bigl)^{p}\Bigr\}^{1/p}\leq A_p\,\|f\|_{L_p(v_{\lambda,\mu})}.
\end{equation}

$(b)$ If $2\leq q<\infty$ and $\phi$ is a function on non-negative integers satisfying
\begin{equation}\label{assumption_for_second_part_of_Hardy-Littlewood_inequality}
\sum\limits_{n=0}^\infty\Bigl((n+1)^{\left(\frac{1}{q'}-\frac{1}{q}\right)\left(\max(\lambda,\mu)+1\right)}|\phi(n)|\Bigl)^{q}<\infty,
\end{equation}
then the algebraic polynomials
\begin{equation*}
\Phi_N(t)=\sum\limits_{n=0}^N \phi(n)\,\widetilde{C}_n^{(\lambda,\mu)}(t)
\end{equation*}
converge in $L_q(v_{\lambda,\mu})$ to a function $f$ satisfying $\hat{f}_n=\phi(n)$, $n=0,1,\ldots$, and
\begin{equation}\label{second_part_of_Hardy-Littlewood_inequality}
\|f\|_{L_q(v_{\lambda,\mu})}\leq A_{q'} \Bigl\{\sum\limits_{n=0}^\infty\Bigl((n+1)^{\left(\frac{1}{q'}-\frac{1}{q}\right)\left(\max(\lambda,\mu)+1\right)}|\phi(n)|\Bigl)^{q}\Bigr\}^{1/q}.
\end{equation}
\end{teoen}

\proofen Let $\sigma=\max(\lambda,\mu)+1$.

(a)  To prove \eqref{first_part_of_Hardy-Littlewood_inequality}, we note that for $p=2$ the Parseval identity implies equality in \eqref{first_part_of_Hardy-Littlewood_inequality} with $A_2=1$. Consider \eqref{first_part_of_Hardy-Littlewood_inequality} as the transformation from $L_p(v_{\lambda,\mu})$ into the sequence $\bigl\{(n+1)^{\sigma}\hat{f}_n\bigr\}_{n=0}^\infty$ in the $\ell_p$ norm with the weight $\bigl\{(n+1)^{-2\sigma}\bigr\}_{n=0}^\infty$ and show that this transformation is of weak type $(1,1)$. We have
\begin{equation*}
m\bigl\{n\colon\, (n+1)^{\sigma}|\hat{f}_n|>t\bigr\}=\sum\limits_{(n+1)^{\sigma}|\hat{f}_n|>t} (n+1)^{-2\sigma}\equiv I_t.
\end{equation*}
By Theorem \ref{main_result_theorem_of_max_value_preprint_for_orthonormal_system}, $|\hat{f}_n|\leq C_1 \|f\|_{L_1(v_{\lambda,\mu})} (n+1)^{\sigma-1}$ and consequently
\begin{equation*}
I_t\leq\sum\limits_{(n+1)>A} (n+1)^{-2\sigma},\quad A=C_2 \left(\dfrac{t}{\|f\|_{L_1(v_{\lambda,\mu})}}\right)^{\frac{1}{2\sigma-1}}.
\end{equation*}
Hence, using the easily verified inequality
\begin{equation*}
\sum\limits_{(n+1)>\widetilde{A}} (1+n)^{-\delta}\leq 2^{\delta-1}\,{\widetilde{A}}^{-\delta+1},\quad \widetilde{A}>0,\quad \delta\geq 2,
\end{equation*}
we observe that, for $\widetilde{A}=A$ and $\delta=2\sigma$,
\begin{equation*}
I_t\leq C_3 \frac{\|f\|_{L_1(v_{\lambda,\mu})}}{t}.
\end{equation*}

The last estimate is a weak $(1,1)$ estimate which, using the Marcinkiewicz interpolation theorem, implies \eqref{first_part_of_Hardy-Littlewood_inequality}.

(b) We have $1<q'\leq2$. For brevity, write $\psi_n$ in place of
$
\Bigl((n+1)^{\left(\frac{1}{q'}-\frac{1}{q}\right)\sigma}|\phi(n)|\Bigl)^{q}.
$ Suppose that $g\in L_{q'}(v_{\lambda,\mu})$ and that $N<N'$ are positive integers. Applying H\"{o}lder's inequality and (a), we find that
\begin{equation}\label{first_inequality_for_second_part_of_Hardy-Littlewood_inequality}
\begin{split}
\Bigl|\int\nolimits_{-1}^1 \Phi_N(t)\, g(t)\, v_{\lambda,\mu}(t)\,dt\Bigr|&=\Bigl|\sum\limits_{n=0}^N\phi(n)\hat{g}_n\Bigr|=\\
&=\Bigl|\sum\limits_{n=0}^N (n+1)^{\left(\frac{1}{q'}-\frac{1}{q}\right)\sigma}\phi(n) \, (n+1)^{\left(\frac{1}{q}-\frac{1}{q'}\right)\sigma}\hat{g}_n\Bigr|\leq\\
&\leq \Bigl\{\sum\limits_{n=0}^N \psi_n\Bigr\}^{1/q}\,\Bigl\{\sum\limits_{n=0}^N \Bigl((n+1)^{\left(\frac{1}{q}-\frac{1}{q'}\right)\sigma}|\hat{g}_n|\Bigl)^{q'} \Bigr\}^{1/q'}\leq\\
&\leq \Bigl\{\sum\limits_{n=0}^N \psi_n\Bigr\}^{1/q} \, A_{q'} \|g\|_{L_{q'}(v_{\lambda,\mu})}.
\end{split}
\end{equation}
Similarly,
\begin{equation}\label{second_inequality_for_second_part_of_Hardy-Littlewood_inequality}
\Bigl|\int\nolimits_{-1}^1 \bigl(\Phi_N(t)-\Phi_{N'}(t)\bigr)\, g(t)\, v_{\lambda,\mu}(t)\,dt\Bigr|\leq \Bigl\{\sum\limits_{n=N+1}^{N'}\psi_n\Bigr\}^{1/q}\, A_{q'} \|g\|_{L_{q'}(v_{\lambda,\mu})}.
\end{equation}
Hence, by \cite[Theorem~(12.13)]{hewitt_ross_book_analysis_1963}, the inequalities \eqref{first_inequality_for_second_part_of_Hardy-Littlewood_inequality} and \eqref{second_inequality_for_second_part_of_Hardy-Littlewood_inequality} lead respectively to the estimates
\begin{equation}\label{third_inequality_for_second_part_of_Hardy-Littlewood_inequality}
\|\Phi_N\|_{L_q(v_{\lambda,\mu})}\leq \Bigl\{\sum\limits_{n=0}^N \psi_n\Bigr\}^{1/q} \, A_{q'}
\end{equation}
and
\begin{equation*}\label{fourth_inequality_for_second_part_of_Hardy-Littlewood_inequality}
\|\Phi_N-\Phi_{N'}\|_{L_q(v_{\lambda,\mu})}\leq \Bigl\{\sum\limits_{n=N+1}^{N'}\psi_n\Bigr\}^{1/q}\, A_{q'}.
\end{equation*}
The last inequality combined with \eqref{assumption_for_second_part_of_Hardy-Littlewood_inequality} show that the sequence $\{\Phi_N\}_{N=1}^\infty$ is a Cauchy sequence in $L_q(v_{\lambda,\mu})$ and therefore convergent in $L_q(v_{\lambda,\mu})$; let $f$ be its limit. Then, by mean convergence,
\begin{equation*}
\hat{f}_n=\lim\limits_{N\to\infty} \bigl(\widehat{\Phi_N}\bigr)_n,\quad n=0,1,\ldots,
\end{equation*}
which is easily seen to equal $\phi(n)$. Moreover, the defining relation
\begin{equation*}
f=\lim\limits_{N\to\infty} \Phi_N\qquad \text{in \, $L_q(v_{\lambda,\mu})$}
\end{equation*}
and the inequality \eqref{third_inequality_for_second_part_of_Hardy-Littlewood_inequality} show that \eqref{second_part_of_Hardy-Littlewood_inequality} holds and so complete the proof.
\hfill$\square$

\section{Hausdorff\,--\,Young-type inequalities\\ for generalized Gegenbauer expansions}\label{section_for_Hausdorff-Young_inequality}

To prove the following result, we use the Riesz\,--\,Thorin interpolation theorem.

\begin{teoen}\label{Hausdorff-Young_inequality_for_generalized_Gegenbauer_expansions}
$(a)$ If $1<p\leq2$ and $f\in L_{p}(v_{\lambda,\mu})$, then
\begin{equation}\label{first_part_of_Hausdorff-Young_inequality}
\Bigl\{\sum\limits_{n=0}^\infty\Bigl((n+1)^{\left(\frac{1}{p'}-\frac{1}{p}\right)\max(\lambda,\mu)}|\hat{f}_n|\Bigr)^{p'}\Bigr\}^{1/p'}\leq B_p\,\|f\|_{L_p(v_{\lambda,\mu})}.
\end{equation}

$(b)$ If $2\leq q<\infty$ and $\phi$ is a function on non-negative integers satisfying
\begin{equation}\label{assumption_for_second_part_of_Hausdorff-Young_inequality}
\sum\limits_{n=0}^\infty\Bigl((n+1)^{\left(\frac{1}{q'}-\frac{1}{q}\right)\max(\lambda,\mu)}|\phi(n)|\Bigl)^{q'}<\infty,
\end{equation}
then the algebraic polynomials
\begin{equation*}
\Phi_N(t)=\sum\limits_{n=0}^N \phi(n)\,\widetilde{C}_n^{(\lambda,\mu)}(t)
\end{equation*}
converge in $L_q(v_{\lambda,\mu})$ to a function $f$ satisfying $\hat{f}_n=\phi(n)$, $n=0,1,\ldots$, and
\begin{equation}\label{second_part_of_Hausdorff-Young_inequality}
\|f\|_{L_q(v_{\lambda,\mu})}\leq B_{q'} \Bigl\{\sum\limits_{n=0}^\infty\Bigl((n+1)^{\left(\frac{1}{q'}-\frac{1}{q}\right)\max(\lambda,\mu)}|\phi(n)|\Bigl)^{q'}\Bigr\}^{1/q'}.
\end{equation}
\end{teoen}

\proofen Let $\sigma=\max(\lambda,\mu)$.

(a) Note that for $p=2$ the Parseval identity implies equality in \eqref{first_part_of_Hausdorff-Young_inequality} with $B_2=1$. We now consider \eqref{first_part_of_Hausdorff-Young_inequality} as the transformation from $L_p(v_{\lambda,\mu})$ into the sequence $\bigl\{(n+1)^{-\sigma}\hat{f}_n\bigr\}_{n=0}^\infty$ in the $\ell_p$ norm with the weight $\bigl\{(n+1)^{2\sigma}\bigr\}_{n=0}^\infty$ and show that this transformation is of strong type $(1,\infty)$. Using Theorem \ref{main_result_theorem_of_max_value_preprint_for_orthonormal_system}, we get
\begin{equation*}
\sup\limits_{n=0,1,\ldots} \Bigl\{(n+1)^{-\sigma} |\hat{f}_n|\Bigr\}\leq B_1 \|f\|_{L_1(v_{\lambda,\mu})}.
\end{equation*}
Thus, applying the Riesz\,--\,Thorin theorem, we deduce \eqref{first_part_of_Hausdorff-Young_inequality} with $B_p=B_1^{\frac{1}{p}-\frac{1}{p'}}$.

(b) The proof of this part is closely related to the proof of part (b) of Theorem \ref{Hardy-Littlewood_inequality_for_generalized_Gegenbauer_expansions}. One can obtain this proof. We give it here for completeness.

We have $1<q'\leq 2$. For brevity, write $\psi_n$ in place of $\Bigl((n+1)^{\left(\frac{1}{q'}-\frac{1}{q}\right)\sigma}|\phi(n)|\Bigl)^{q'}$. Suppose that $g\in L_{q'}(v_{\lambda,\mu})$ and that $N<N'$ are positive integers. Applying H\"{o}lder's inequality and (a), we find that
\begin{equation}\label{first_inequality_for_second_part_of_Hausdorff-Young_inequality}
\begin{split}
\Bigl|\int\nolimits_{-1}^1 \Phi_N(t)\, g(t)\, v_{\lambda,\mu}(t)\,dt\Bigr|&=\Bigl|\sum\limits_{n=0}^N\phi(n)\hat{g}_n\Bigr|=\\
&=\Bigl|\sum\limits_{n=0}^N (n+1)^{\left(\frac{1}{q'}-\frac{1}{q}\right)\sigma}\phi(n) \, (n+1)^{\left(\frac{1}{q}-\frac{1}{q'}\right)\sigma}\hat{g}_n\Bigr|\leq\\
&\leq \Bigl\{\sum\limits_{n=0}^N \psi_n\Bigr\}^{1/q'}\,\Bigl\{\sum\limits_{n=0}^N \Bigl((n+1)^{\left(\frac{1}{q}-\frac{1}{q'}\right)\sigma}|\hat{g}_n|\Bigl)^{q} \Bigr\}^{1/q}\leq\\
&\leq \Bigl\{\sum\limits_{n=0}^N \psi_n\Bigr\}^{1/q'} \, B_{q'} \|g\|_{L_{q'}(v_{\lambda,\mu})}.
\end{split}
\end{equation}
Similarly,
\begin{equation}\label{second_inequality_for_second_part_of_Hausdorff-Young_inequality}
\Bigl|\int\nolimits_{-1}^1 \bigl(\Phi_N(t)-\Phi_{N'}(t)\bigr)\, g(t)\, v_{\lambda,\mu}(t)\,dt\Bigr|\leq \Bigl\{\sum\limits_{n=N+1}^{N'}\psi_n\Bigr\}^{1/q'}\, B_{q'} \|g\|_{L_{q'}(v_{\lambda,\mu})}.
\end{equation}
Hence, by \cite[Theorem~(12.13)]{hewitt_ross_book_analysis_1963}, the inequalities \eqref{first_inequality_for_second_part_of_Hausdorff-Young_inequality} and \eqref{second_inequality_for_second_part_of_Hausdorff-Young_inequality} lead respectively to the estimates
\begin{equation}\label{third_inequality_for_second_part_of_Hausdorff-Young_inequality}
\|\Phi_N\|_{L_q(v_{\lambda,\mu})}\leq \Bigl\{\sum\limits_{n=0}^N \psi_n\Bigr\}^{1/q'} \, B_{q'}
\end{equation}
and
\begin{equation*}\label{fourth_inequality_for_second_part_of_Hausdorff-Young_inequality}
\|\Phi_N-\Phi_{N'}\|_{L_q(v_{\lambda,\mu})}\leq \Bigl\{\sum\limits_{n=N+1}^{N'}\psi_n\Bigr\}^{1/q'}\, B_{q'}.
\end{equation*}
The last inequality combined with \eqref{assumption_for_second_part_of_Hausdorff-Young_inequality} show that the sequence $\{\Phi_N\}_{N=1}^\infty$ is a Cauchy sequence in $L_q(v_{\lambda,\mu})$ and therefore convergent in $L_q(v_{\lambda,\mu})$; let $f$ be its limit. Then, by mean convergence,
\begin{equation*}
\hat{f}_n=\lim\limits_{N\to\infty} \bigl(\widehat{\Phi_N}\bigr)_n,\quad n=0,1,\ldots,
\end{equation*}
which is easily seen to equal $\phi(n)$. Moreover, the defining relation
\begin{equation*}
f=\lim\limits_{N\to\infty} \Phi_N\qquad \text{in \, $L_q(v_{\lambda,\mu})$}
\end{equation*}
and the inequality \eqref{third_inequality_for_second_part_of_Hausdorff-Young_inequality} show that \eqref{second_part_of_Hausdorff-Young_inequality} holds and so complete the proof.
\hfill$\square$

\section{Unification of the Hardy\,--\,Littlewood-type\\ and the Hausdorff\,--\,Young-type inequalities}\label{section_for_unification_of_inequalities}

Theorem \ref{unification_of_different_types_of_inequalities} contains the Hardy\,--\,Littlewood-type and the Hausdorff\,--\,Young-type inequalities for the expansions by orthonormal polynomials with respect to the weight function $v_{\lambda,\mu}$ (see \eqref{weight_function}). To prove it, we need Stein's modification of the Riesz\,--\,Thorin interpolation theorem (see \cite[Theorem~2, p.~485]{stein_article_interpolation_1956}) given below.

\begin{teoen}[(Stein)]\label{Stein's_modification}
Suppose $\nu_1$ and $\nu_2$ are $\sigma$-finite measures on $M$ and $S$, respectively, and $T$ is a linear operator defined on $\nu_1$-measurable functions on $M$ to $\nu_2$-measurable functions on $S$. Let $1\leq r_0,\,r_1,\,s_0,\,s_1\leq\infty$ and $\frac{1}{r}=\frac{1-t}{r_0}+\frac{t}{r_1}$, $\frac{1}{s}=\frac{1-t}{s_0}+\frac{t}{s_1}$, where $0\leq t\leq1$. Suppose further that
\begin{equation*}
\|(Tg) \cdot v_i\|_{L_{s_i}(S,\nu_2)}\leq L_i\|g \cdot u_i\|_{L_{r_i}(M,\nu_1)},\quad i=0,1,
\end{equation*}
where $u_i$ and $v_i$ are non-negative weight functions. Let $u=u_0^{1-t} \cdot u_1^t$, $v=v_0^{1-t} \cdot v_1^t$.

Then
\begin{equation*}
\|(Tg) \cdot v\|_{L_{s}(S,\nu_2)}\leq L\|g \cdot u\|_{L_{r}(M,\nu_1)}
\end{equation*}
with $L=L_0^{1-t} \cdot L_1^t$.
\end{teoen}

\begin{teoen}\label{unification_of_different_types_of_inequalities}
Let $\sigma=\max(\lambda,\mu)$.

$(a)$ If $1<p\leq 2$, $f\in L_p(v_{\lambda,\mu})$, and $p\leq s\leq p'$, then
\begin{equation}\label{first_part_of_unification}
\Bigl\{\sum\limits_{n=0}^\infty\Bigl((n+1)^{\left(\frac{1}{s}-\frac{1}{p}\right)\sigma+\left(\frac{1}{p'}-\frac{1}{s}\right)(\sigma+1)}|\hat{f}_n|\Bigr)^s\Bigr\}^{1/s}\leq C_p(s)\,\|f\|_{L_p(v_{\lambda,\mu})}.
\end{equation}

$(b)$ If $2\leq q<\infty$, $q'\leq r\leq q$, and $\phi$ is a function on non-negative integers satisfying
\begin{equation*}
\sum\limits_{n=0}^\infty\Bigl((n+1)^{\left(\frac{1}{q'}-\frac{1}{r}\right)\sigma+\left(\frac{1}{r}-\frac{1}{q}\right)(\sigma+1)}|\phi(n)|\Bigr)^{r'}<\infty,
\end{equation*}
then the algebraic polynomials
\begin{equation*}
\Phi_N(t)=\sum\limits_{n=0}^N \phi(n)\widetilde{C}_n^{(\lambda,\mu)}(t)
\end{equation*}
converge in $L_q(v_{\lambda,\mu})$ to a function $f$ satisfying $\hat{f}_n=\phi(n)$, $n=0,1,\ldots$, and
\begin{equation*}
\|f\|_{L_q(v_{\lambda,\mu})}\leq C_{q'}(r) \Bigl\{\sum\limits_{n=0}^\infty\Bigl((n+1)^{\left(\frac{1}{q'}-\frac{1}{r}\right)\sigma+\left(\frac{1}{r}-\frac{1}{q}\right)(\sigma+1)}|\phi(n)|\Bigr)^{r'}\Bigr\}^{1/r'}.
\end{equation*}

\end{teoen}

\proofen (a) This part was proved for $s=p$ (with $C_p(p)=A_p$) and $s=p'$ (with $C_p(p')=B_p$) in Theorems \ref{Hardy-Littlewood_inequality_for_generalized_Gegenbauer_expansions} and \ref{Hausdorff-Young_inequality_for_generalized_Gegenbauer_expansions}, respectively. So for $p=2$, we obtain the equality in \eqref{first_part_of_unification} with $C_2(2)=1$.

Consider now the case that $1<p<2$. To prove \eqref{first_part_of_unification}, we set in Theorem \ref{Stein's_modification}: $M=[-1,1]$, $\nu_1$ the Lebesgue measure, $S=\{n\}_{n=0}^\infty$, $\nu_2$ the counting measure, $g=f$, $Tg=\{\hat{f}_n\}_{n=0}^\infty$, $r=r_0=r_1=p$, $u=u_0=u_1=v_{\lambda,\mu}$, $s_0=p'$, $s_1=p$, $v_0=\bigl\{(n+1)^{\left(\frac{1}{p'}-\frac{1}{p}\right)\sigma}\bigr\}_{n=0}^\infty$, $v_1=\bigl\{(n+1)^{\left(\frac{1}{p'}-\frac{1}{p}\right)(\sigma+1)}\bigr\}_{n=0}^\infty$, $L_0=B_p$, $L_1=A_p$, and $\frac{1}{s}=\frac{1-t}{p'}+\frac{t}{p}$. As $\frac{1}{p}+\frac{1}{p'}=1$, $\frac{1}{s}-\frac{1}{p}=(1-t)\bigl(\frac{1}{p'}-\frac{1}{p}\bigr)$, $\frac{1}{p'}-\frac{1}{s}=t\bigl(\frac{1}{p'}-\frac{1}{p}\bigr)$, the proof of \eqref{first_part_of_unification} is concluded.

Because of
\begin{equation*}
1-t=\frac{\frac{1}{s}-\frac{1}{p}}{\frac{1}{p'}-\frac{1}{p}},\quad t=\frac{\frac{1}{p'}-\frac{1}{s}}{\frac{1}{p'}-\frac{1}{p}},
\end{equation*}
it is clear that $C_p(s)=B_p^{1-t}A_p^t$.

(b) Taking into account the previously given proofs (see parts (b) and (b) in Theorems \ref{Hardy-Littlewood_inequality_for_generalized_Gegenbauer_expansions} and \ref{Hausdorff-Young_inequality_for_generalized_Gegenbauer_expansions}, respectively), the proof is obvious and left to the reader.
\hfill$\square$

\section*{Acknowledgements}

This work was done thanks to the remarkable papers \cite{ditzian_article_smoothness_2012,ditzian_article_norm_and_smoothness_2015,ditzian_article_estimates_2011} of Z.~Ditzian.

\begin{Biblioen}

\bibitem{andrews_askey_roy_book_special_functions_1999}G.\,E.~Andrews, R.~Askey, and R.~Roy, \textit{Special Functions}, Encyclopedia of Mathematics and its Applications \textbf{71}, Cambridge University Press, Cambridge, 1999.

\bibitem{dai_xu_book_approximation_theory_2013}F.~Dai and Y.~Xu, \textit{Approximation theory and harmonic analysis on spheres and balls}, Springer Monographs in Mathematics, Springer, 2013.

\bibitem{ditzian_article_estimates_2011}Z.~Ditzian, Estimates of the coefficients of the Jacobi expansion by measures of smoothness, \textit{J. Math. Anal. Appl.} \textbf{384} (2011), 303--306.

\bibitem{ditzian_article_smoothness_2012}Z.~Ditzian, Relating smoothness to expressions involving Fourier coefficients or to a Fourier transform, \textit{J. Approx. Theory} \textbf{164} (2012), 1369--1389.

\bibitem{ditzian_article_norm_and_smoothness_2015}Z.~Ditzian, Norm and smoothness of a function related to the coefficients of its expansion, \textit{J. Approx. Theory} \textbf{196} (2015), 101--110.

\bibitem{dunkl_xu_book_orthogonal_polynomials_2014}C.\,F.~Dunkl and Y.~Xu, \textit{Orthogonal polynomials of several variables}, 2nd ed., Encyclopedia of Mathematics and its Applications \textbf{155}, Cambridge University Press, Cambridge, 2014.

\bibitem{hewitt_ross_book_analysis_1963}E.~Hewitt and K.\,A.~Ross, \textit{Abstract harmonic analysis. Vol. I}, Springer-Verlag, Heidelberg, 1963.

\bibitem{stein_article_interpolation_1956}E.~Stein, Interpolation of linear operators, \textit{Trans. Amer. Math. Soc.} \textbf{83} (1956), 482--492.

\bibitem{stein_weiss_article_interpolation_1958}E.~Stein and G.~Weiss, Interpolation of operators with change of measures, \textit{Trans. Amer. Math. Soc.} \textbf{87} (1958), 159--172.

\bibitem{szego_book_orthogonal_polynomials_1975}G.~Szeg\"{o}, \textit{Orthogonal polynomials}, 4th ed., American Mathematical Society Colloquium Publications \textbf{23}, American Mathematical Society, Providence, Rhode Island, 1975.

\bibitem{veprintsev_preprint_max_value_2015}R.\,A.~Veprintsev, On the asymptotic behavior of the maximum absolute value of generalized Gegenbauer polynomials, arXiv preprint 1602.01023 (2015).

\end{Biblioen}

\noindent \textsc{Department of scientific research, Tula State University, Tula, Russia }

\noindent \textit{E-mail address}: \textbf{veprintsevroma@gmail.com}

\end{document}